\documentclass[12pt]{amsart}

\newtheorem{theorem}{Theorem}[section]
\newtheorem{lemma}[theorem]{Lemma}
\newtheorem{corollary}[theorem]{Corollary}

\numberwithin{equation}{section}

\begin{document}

\baselineskip=15.5pt

\title[Principal bundles on a projective variety]{On
principal bundles over a projective variety
defined over a finite field}

\author[I. Biswas]{Indranil Biswas}

\address{School of Mathematics, Tata Institute of Fundamental
Research, Homi Bhabha Road, Bombay 400005, India}

\email{indranil@math.tifr.res.in}



\date{}

\begin{abstract}
Let $M$ be a geometrically irreducible smooth projective
variety, defined over a finite field $k$, such that
$M$ admits a $k$--rational point $x_0$. Let $\varpi(M,x_0)$
denote the corresponding fundamental group--scheme
introduced by Nori. Let $E_G$ be
a principal $G$--bundle over $M$, where $G$ is a reduced
reductive linear algebraic group defined over the field $k$.
Fix a polarization $\xi$ on $M$. We
prove that the following three statements are equivalent:
\begin{enumerate}
\item The principal $G$--bundle $E_G$ over $M$ is given by a
homomorphism $\varpi(M,x_0)\, \longrightarrow\,G$.

\item There are integers $b\, >\, a\, \geq\, 1$, such that
the principal $G$--bundle $(F^b_M)^*E_G$ is isomorphic to
$(F^a_M)^*E_G$, where $F_M$ is the absolute Frobenius morphism
of $M$.

\item  The principal $G$--bundle $E_G$ is 
strongly semistable,
$\text{degree}(c_2(\text{ad}(E_G))c_1(\xi)^{d-2})\, =\,0$,
where $d\, =\, \dim M$,
and $\text{degree}(c_1(E_G(\chi))c_1(\xi)^{d-1})\,=\, 0$
for every character $\chi$ of $G$, where
$E_G(\chi)$ is the line bundle over $M$ associated to $E_G$
for $\chi$.
\end{enumerate}
In \cite{Su}, the equivalence between the first statement and the
third statement was proved under the extra
assumption that $\dim M\,=\,1$ and $G$ is semisimple.

\end{abstract}

\maketitle

\section{Introduction}\label{sec1}

Let $M$ be a geometrically irreducible smooth projective
variety, defined over a finite field $k$, such that $M$ admits
a $k$--rational point. Fix a $k$--rational point $x_0$ of $M$.
Also, fix a polarization $\xi$ on $M$ in order to
define the degree of torsionfree coherent sheaves on $M$.

Let $G$ be a reduced reductive linear algebraic group defined over
the field $k$.
We recall that a principal $G$--bundle $E_G$ over
$M$ is defined to be \textit{strongly semistable} if $(F^n_M)^*E_G$
is semistable for all $n\, \geq\,1$, where
$$
F_M\, :\, M\, \longrightarrow\, M
$$
is the absolute Frobenius
morphism (see \cite[p. 129, Definition 1.1]{Ra} and
\cite[pp. 131--132, Lemma 2.1]{Ra} for the definition
of a semistable principal $G$--bundle). It is known that
a principal $G$--bundle $E_G$ is strongly semistable if
and only if for every representation
$$
\eta\, :\, G\, \longrightarrow\, \text{GL}(V)\, ,
$$
where $V$ is a finite dimensional $k$--vector
space, satisfying the condition that $\eta$
takes the connected component of the center of $G$ to
the center of  $\text{GL}(V)$, the vector bundle
$$
E_G(V)\, :=\, E\times^G V\, \longrightarrow\, M
$$
associated to $E_G$ for $\eta$
is semistable (see \cite[p. 288, Theorem 3.23]{RR}).
We are using the notation
$E_G\times^G V$ of \cite[p. 114, D\'efinition 1.3.1]{Gi}
in order to denote the
quotient of $E_G\times V$ by the ``twisted'' diagonal
action of $G$; this notation will be used throughout. 

The notion of fundamental group--scheme
was introduced by Nori \cite{No0}.
The fundamental
group--scheme of $M$ with base point $x_0$ will be denoted
by $\varpi(M,x_0)$. We recall that $\varpi(M,x_0)$ is
the affine group--scheme over $k$ associated to
the neutral Tannakian category defined by the
essentially finite vector bundles on $M$; the fiber
functor of the neutral Tannakian category
sends an essentially finite vector bundle to its
fiber over $x_0$ (see Section \ref{sec3}).

There is a tautological principal $\varpi(M,x_0)$--bundle
$E_{\varpi(M,x_0)}$ over $M$; its description is recalled
in Section \ref{sec3}. Therefore, given any homomorphism
\begin{equation}\label{e0}
\rho\, :\, \varpi(M,x_0)\, \longrightarrow\, G
\end{equation}
we have the principal $G$--bundle
\begin{equation}\label{eb}
E_\rho\, :=\, E_{\varpi(M,x_0)}
\times^{\varpi(M,x_0)} G \,\longrightarrow\, M
\end{equation}
which is obtained by extending
the structure group of the tautological
principal $\varpi(M,x_0)$--bundle
$E_{\varpi(M,x_0)}$ using $\rho$.

For any zero--cycle $C\,\in\, \text{CH}^0(M)$, let $[C]\, \in\,
\mathbb Z$ be the degree of $C$.

Our aim is to prove the following theorem.

\begin{theorem}\label{thm0}
Let $E_G$ be a principal $G$--bundle over $M$.
The following three statements are equivalent.
\begin{enumerate}
\item There are integers $b\, >\, a\, \geq\, 1$, such that
the principal $G$--bundle $(F^b_M)^*E_G$ is isomorphic to
$(F^a_M)^*E_G$.

\item There is a homomorphism $\rho$ as in Eq. \eqref{e0} such
that $E_G$ is isomorphic to the principal $G$--bundle
$E_\rho$ defined in Eq. \eqref{eb}.

\item The principal $G$--bundle $E_G$ is strongly semistable
and the following two conditions hold:
\begin{itemize}
\item for each character $\chi$ of $G$, the line bundle
$E_G(\chi)\, =\, E_G\times^G k\, \longrightarrow\, M$
associated to $E_G$ for $\chi$ has the property that
$$
[c_1(E_G(\chi))c_1(\xi)^{d-1}]\, \in\, {\mathbb Z}
$$
vanishes, where $d\,=\, \dim M$, and

\item $[c_2({\rm ad}(E_G))c_1(\xi)^{d-2}]\, =\, 0$, where
${\rm ad}(E_G)$ is the adjoint vector bundle of $E_G$.
\end{itemize}
\end{enumerate}
\end{theorem}

In \cite{Su}, the equivalence between the second and the
third statements in Theorem \ref{thm0} was
proved under the extra assumption that $\dim M\,=\,1$ and
$G$ is semisimple.

We also prove the following:

\begin{corollary}\label{cor0}
A vector bundle $E$ over $M$ is essentially finite if and only
if there are integers $b\, >\, a\, \geq\, 1$, such that
the vector bundle $(F^b_M)^*E$ is isomorphic to $(F^a_M)^*E$.
\end{corollary}

\section{Strongly semistable principal bundles}\label{sec2}

Let $k$ be a finite field.
Let $M$ be a geometrically irreducible smooth projective
variety defined over $k$. Let $d$ be the dimension of $M$.

Fix a very ample line bundle $\xi$ over $M$. For any
torsionfree coherent sheaf $V$ on $M$, define
$$
\text{degree}(V)\, :=\, [c_1(V)c_1(\xi)^{d-1}]\, \in\,
\mathbb Z\, .
$$

A vector bundle $E$ over $M$ is called \textit{finite}
if there are two distinct polynomials
$$
f\, , g\, \in\, {\mathbb Z}[X]
$$
with nonnegative coefficients such that the vector
bundle $f(E)$ is isomorphic to $g(E)$
(see \cite[p. 80, Lemma 3.1]{No}). A vector bundle
$V$ over $M$ of degree zero
is called \textit{essentially finite}
if there is a finite vector bundle $E$ and a subbundle
$W\, \subset\, E$ such that
\begin{itemize}
\item the vector bundle $W$ is of degree zero,
and
\item $V$ is isomorphic to a quotient of $W$.
\end{itemize}
(See \cite[p. 82, Definition]{No}.)

Let $G$ be a reduced reductive linear algebraic group
defined over $k$. For a principal $G$--bundle $E_G$ over $M$,
the \textit{adjoint vector bundle} $\text{ad}(E_G)$
is the one associated to $E_G$ for the adjoint action
of $G$ on its own Lie algebra.

\begin{lemma}\label{lem1}
Let $E_G$ be a strongly semistable principal $G$--bundle over
$M$ satisfying the following two conditions:
\begin{itemize}
\item for each character $\chi$ of $G$, the line bundle
$E_G(\chi)\, =\, E_G\times^G k\,\longrightarrow\, M$
associated to $E_G$ for $\chi$ has the property that
$$
[c_1(E_G(\chi))c_1(\xi)^{d-1}]\, \in\, {\mathbb Z}
$$
vanishes, and

\item $[c_2({\rm ad}(E_G))c_1(\xi)^{d-2}]\, \in\,
{\mathbb Z}$ vanishes.
\end{itemize}
Let $V$ be a finite dimensional representation
of $G$ over $k$.
Then the associated vector bundle $E_G(V)\, =\,
E_G\times^G V$ over $M$ is essentially finite.
\end{lemma}

\begin{proof}
Since $E_G$ is strongly semistable, the associated vector bundle
$E_G(V)$ is also strongly semistable \cite[p. 288,
Theorem 3.23]{RR}. Taking the character $\chi$ in the statement
of the lemma to be the one associated to the
representation $\bigwedge^{\text{top}}V$ of $G$ we conclude that
$[c_1(E_G(V))c_1(\xi)^{d-1}]\, =\, 0$.

We will show that $[c_2(E_G(V))c_1(\xi)^{d-2}]\,=\, 0$.
For that, fix a filtration of $G$--modules
\begin{equation}\label{e1}
0\,=\, V_0\, \subset\, V_1\, \subset\, \cdots\, \subset\,
V_{n-1}\, \subset\, V_n\, :=\, V
\end{equation}
such that each successive quotient $V_i/V_{i-1}$, $i\, \in\,
[1\, , n]$, is irreducible.

The center of $G$ will be denoted by $Z$, and the Lie algebra
of $G$ will be denoted by $\mathfrak g$.
We note that $\mathfrak g$ is a faithful $G/Z$--module.

Since $V_i/V_{i-1}$ in Eq. \eqref{e1} is an irreducible $G$--module,
the action of the center $Z$ on $V_i/V_{i-1}$ is multiplication by
a fixed character $\chi_i$ of $G$. Hence $Z$ acts trivially on
$$
\text{End}(V_i/V_{i-1})\, =\, (V_i/V_{i-1})\bigotimes
(V_i/V_{i-1})^*\, .
$$
Therefore, the action of $G$ on $\text{End}(V_i/V_{i-1})$
gives an action of $G/Z$ on $\text{End}(V_i/V_{i-1})$.
We will consider $\text{End}(V_i/V_{i-1})$ as a
$G/Z$--module. Since 
$\mathfrak g$ is a faithful $G/Z$--module, each
$G/Z$--module $V_i/V_{i-1}$ , $i\, \in\,
[1\, , n]$, is a subquotient of a $G/Z$--module of the form
\begin{equation}\label{e2}
{\mathcal V}\, :=\, \bigoplus_{i=1}^m
{\mathfrak g}^{\otimes a_i}\otimes ({\mathfrak g}^*)^{\otimes 
b_i}\, ,
\end{equation}
where $a_i\, ,b_i$ are nonnegative integers and $m\, >\, 0$
(see \cite[Proposition 4.4]{BPS}). We recall that a
subquotient of ${\mathcal V}$ is a sub $G/Z$--module
of a quotient $G/Z$--module of $\mathcal V$ (which is same
as being a quotient of a sub module of $\mathcal V$). 

Using \cite[p. 288, Theorem 3.23]{RR} it follows that
the adjoint bundle $\text{ad}(E_G)$, which is associated to
$E_G$ for the $G$--module $\mathfrak g$, is strongly semistable.
Also, $\bigwedge^{\text{top}}E_G({\mathcal V})$ is a trivial line
bundle because the $G$--module $\bigwedge^{\text{top}}{\mathcal V}$
is trivial. Hence the vector bundle $E_G({\mathcal V})$ associated
to $E_G$ for the $G$--module $\mathcal V$ in Eq. \eqref{e2} has the
following properties:
\begin{itemize}
\item $E_G({\mathcal V})$ is strongly semistable,

\item $c_1(E_G({\mathcal V}))\, =\, 0$, and

\item $[c_2(E_G({\mathcal V}))c_1(\xi)^{d-2}]\, =\, 0$
(recall that $[c_2({\rm ad}(E_G))c_1(\xi)^{d-2}]\,=\, 0$ by assumption).
\end{itemize}
Since $\text{End}(V_i/V_{i-1})$ is a subquotient of $\mathcal V$,
there is a quotient $G/Z$--module
$$
Q\, =\, {\mathcal V}/{\mathcal K}
$$
of $\mathcal V$ such that
the $G/Z$--module $\text{End}(V_i/V_{i-1})$ is isomorphic to
a submodule of $Q$. Let $E_G({\mathcal K})$,
$E_G(\text{End}(V_i/V_{i-1}))$ and 
$E_G(Q/\text{End}(V_i/V_{i-1}))$ be the vector bundles associated
to $E_G$ for the $G$--modules ${\mathcal K}$,
$\text{End}(V_i/V_{i-1})$ and $Q/\text{End}(V_i/V_{i-1})$
respectively. So, as elements of the Grothendieck $K$--group
$K(M)$,
\begin{equation}\label{gk}
E_G({\mathcal V})\, =\, E_G({\mathcal K}) + 
E_G(Q/\text{End}(V_i/V_{i-1})) +E_G(\text{End}(V_i/V_{i-1}))
\, \in\, K(M)\, .
\end{equation}
Since $E_G$ is strongly semistable,
using \cite[p. 288, Theorem 3.23]{RR} we know that the
associated vector
bundles $E_G({\mathcal K})$, $E_G(Q/\text{End}(V_i/V_{i-1}))$ and
$E_G(\text{End}(V_i/V_{i-1}))$ are strongly semistable. Also
$\bigwedge^{\text{top}} E_G({\mathcal K})$, 
$\bigwedge^{\text{top}}E_G(Q/\text{End}(V_i/V_{i-1}))$ and
$\bigwedge^{\text{top}}E_G(\text{End}(V_i/V_{i-1}))$ are trivial
line bundles because ${\mathcal K}$,
$\text{End}(V_i/V_{i-1})$ and $Q/\text{End}(V_i/V_{i-1})$ are all
given by $G/Z$--modules, and $G/Z$ does not have any nontrivial
character. From Bogomolov's inequality,
\cite[p. 500, Theorem]{Bo} (see also
\cite[p. 252, Theorem 0.1]{Lan}), we know that for a
strongly semistable vector bundle $W\, \longrightarrow\, M$
with $\bigwedge^{\text{top}}W$ trivial,
$$
[c_2(W)c_1(\xi)^{d-2}]\, \geq\, 0\, .
$$
Since $[c_2(E_G({\mathcal V}))c_1(\xi)^{d-2}]\, =\, 0$, from
Eq. \eqref{gk} we now have
$$
[c_2(E_G(\text{End}(V_i/V_{i-1})))c_1(\xi)^{d-2}]\, \geq\, 0\, .
$$
Consequently, $[c_2(E_G(V))c_1(\xi)^{d-2}]\, 
\in\,{\mathbb Z}$ vanishes.

The semistable vector bundle $E\, \longrightarrow\, M$
of fixed rank with
$$
[c_1(E)c_1(\xi)^{d-1}]\, =\, 0\, =\,
[c_2(E)c_1(\xi)^{d-2}]
$$
form a bounded family \cite[p. 269, Theorem 4.2]{Lan}.
The field $k$ being finite, it follows from the above
boundedness theorem that the set of
vector bundles $\{(F^n_M)^*E_G(V)\}_{n\geq 1}$ contains only
finitely many isomorphism classes. Hence there are
integers $b\, >\, a\, \geq\, 1$ such that
\begin{equation}\label{e3}
(F^b_M)^*E_G(V)\, =\, (F^a_M)^*E_G(V)\, .
\end{equation}

Since $(F^{b-a}_M)^*(F^a_M)^*E_G(V)\, =\,
(F^b_M)^*E_G(V)$, from Eq. \eqref{e3} we have
\begin{equation}\label{eqf}
(F^{b-a}_M)^*(F^a_M)^*E_G(V)\, =\, (F^a_M)^*E_G(V)\,  .
\end{equation}
Consequently, from a theorem of Lange--Stuhler and
Deligne we conclude the following:

There is an \'etale Galois covering
$$
\phi\, :\, Y\, \longrightarrow\, M
$$
such that the pull back $\phi^*(F^a_M)^*E_G(V)$ is a
trivial vector bundle; see \cite[p. 75, Theorem 1.4]{LS}
and \cite[\S~3.2]{La}.

Now from \cite[pp. 552--553, Proposition 2.3]{BH}
it follows that the vector bundle
$E_G(V)$ is essentially finite. This
completes the proof of the lemma.
\end{proof}

\section{Fundamental group--scheme and principal
bundles}\label{sec3}

Henceforth, we will assume that
the variety $M$ admits a $k$--rational
point. Fix a $k$--rational point $x_0$ of $M$.

The fundamental
group--scheme of $M$ with base point $x_0$ will be denoted
by $\varpi(M,x_0)$ (see \cite[p. 85, Definition 1 and
Proposition 2]{No}). There is a tautological principal
$\varpi(M,x_0)$--bundle over $M$ whose construction
is recalled below. (See also \cite[p. 84, Definition]{No}
for this tautological principal $\varpi(M,x_0)$--bundle.)

The fundamental group--scheme $\varpi(M,x_0)$ is defined
by giving the corresponding neutral Tannakian category.
More precisely,
consider the neutral Tannakian category defined by the
essentially finite vector bundles over $M$ (their definition
was recalled in Section \ref{sec2}); the fiber functor
for the neutral Tannakian category
sends an essentially finite vector bundle $W$
to the $k$--vector space $W_{x_0}$, where $W_{x_0}$
is the fiber of $W$ over the base point $x_0$. 
The fundamental group--scheme $\varpi(M,x_0)$ is defined
to be the group--scheme
associated to this neutral Tannakian category.
Consequently, each essentially finite vector bundle
$W$ over $M$ gives a finite dimensional
representation of $\varpi(M,x_0)$ over $k$ such that the
underlying $k$--vector space is $W_{x_0}$.

Let $\text{Rep}(\varpi(M,x_0))$ denote the neutral Tannakian
category defined be the representation
of $\varpi(M,x_0)$. So $\text{Rep}(\varpi(M,x_0))$
is equivalent to the above neutral Tannakian category
defined by the essentially finite vector bundles over $M$.
Let $\text{Vect}(M)$ denote the category
of vector bundles over $X$. We have a tautological functor
\begin{equation}\label{e5}
{\mathcal F}\,:\, \text{Rep}(\varpi(M,x_0))\, \longrightarrow\,
\text{Vect}(M)
\end{equation}
that sends any $W\, \in\, \text{Rep}(\varpi(M,x_0))$
to the essentially finite vector bundle $W$.
Using \cite[p. 149, Theorem 3.2]{DM}, \cite[Lemma 2.3, Proposition
2.4]{No}, this functor ${\mathcal F}$ in Eq. \eqref{e5} defines
a principal $\varpi(M,x_0)$--bundle over $M$.

The above tautological
principal $\varpi(M,x_0)$--bundle over $M$ will be denoted
by $E_{\varpi(M,x_0)}$.

We also note that the restriction of the principal 
$\varpi(M,x_0)$--bundle $E_{\varpi(M,x_0)}$ to the
base point $x_0$
is canonically trivialized. This trivialization
is obtained from the facts that
the fiber functor for $\text{Rep}(\varpi(M,x_0))$ takes
any $W$ to the fiber of ${\mathcal F}(W)$ over $x_0$
(see Eq. \eqref{e5} for $\mathcal F$), and
$E_{\varpi(M,x_0)}$ is defined by ${\mathcal F}$. Let
\begin{equation}\label{e4}
e_0\, \in\, (E_{\varpi(M,x_0)})_{x_0}
\end{equation}
be the point that corresponds to the identity element in
$\varpi(M,x_0)$ by the canonical trivialization of
the $\varpi(M,x_0)$--torsor $(E_{\varpi(M,x_0)})_{x_0}$.

As in Section \ref{sec1}, for any homomorphism
\begin{equation}\label{e-1}
\rho\, :\, \varpi(M,x_0)\, \longrightarrow\, G
\end{equation}
of group--schemes,
the principal $G$--bundle over $M$ obtained by extending the
structure group of the principal
$\varpi(M,x_0)$--bundle $E_{\varpi(M,x_0)}$ using $\rho$
will be denoted by $E_\rho$.

\begin{theorem}\label{thm1}
Let $E_G$ be a principal $G$--bundle over $M$.
The following three statements are equivalent.
\begin{enumerate}
\item There are integers $b\, >\, a\, \geq\, 1$, such that
the principal $G$--bundle $(F^b_M)^*E_G$ is isomorphic to
$(F^a_M)^*E_G$.

\item There is a homomorphism $\rho$ as in Eq. \eqref{e-1} such that
$E_G$ is isomorphic to the principal $G$--bundle $E_\rho$.

\item The principal $G$--bundle $E_G$ is strongly semistable
and the following two conditions hold:
\begin{itemize}
\item for each character $\chi$ of $G$, the line bundle
$E_G(\chi)\, =\, E_G\times^G k
\,\longrightarrow\, M$ associated to $E_G$ for
$\chi$ has the property that
$$
[c_1(E_G(\chi))c_1(\xi)^{d-1}]\, \in\, {\mathbb Z}
$$
vanishes, and

\item $[c_2({\rm ad}(E_G))c_1(\xi)^{d-2}]\, \in\,
{\mathbb Z}$ vanishes, where ${\rm ad}(E_G)$ is the adjoint
bundle of $E_G$.
\end{itemize}
\end{enumerate}
\end{theorem}

\begin{proof}
We will show that
the first statement implies the third statement. Assume that
$(F^b_M)^*E_G$ is isomorphic to $(F^a_M)^*E_G$, where
$b\, >\, a\, \geq\,1$. Then
all the numerical invariants of $E_G$ vanish. To prove that
$E_G$ is strongly semistable, first note that the vector
bundle $(F^b_M)^*\text{ad}(E_G)$ is isomorphic to
$(F^a_M)^*\text{ad}(E_G)$. For a coherent subsheaf
$E\, \subset\, \text{ad}(E_G)$ with $\text{degree}(E)\, >\, 0$
(the degree is defined using $\xi$), we have
$$
\text{degree}((F^{jb-ja+1}_M)^*E)\, =\,
p(jb-ja+1)\cdot \text{degree}(E)
$$
for the subsheaf
$$
(F^{jb-ja+1}_M)^*E\, \subset\,
(F^{j(b-a)}_M)^* (F^a_M)^*\text{ad}(E_G)\, =\,
(F^a_M)^*\text{ad}(E_G)\, ,
$$
where $p$ is the characteristic of the field $k$ and
$j$ is any positive integer. But any given
vector bundle can not contain subsheaves of arbitrarily large
degrees. In particular, $(F^a_M)^*\text{ad}(E_G)$
does not contain subsheaves of arbitrarily large degrees.
Therefore, we conclude that $\text{ad}(E_G)$
does not contain any subsheaf of positive degree.
Hence $\text{ad}(E_G)$ is semistable. This immediately
implies that the principal $G$--bundle $E_G$ is semistable.
Now replacing $E_G$ by $(F^n_M)^*E_G$ in the above argument
it follows that $E_G$ is strongly semistable.

Therefore, the first statement in the theorem
implies the third statement. We will now show that the
third statement implies the first statement.

The family of principal $G$--bundles
over $M$ satisfying all the
conditions in the third statement in the theorem is
bounded \cite[p. 533, Theorem 7.3]{Lan2} (see
also \cite{HK}). Note that if a principal $G$--bundle
$E'_G$ satisfies all the conditions in the third statement,
then $F^*_M E'_G$ also satisfies all these conditions.
Therefore, from the finiteness of
the field $k$ we conclude that if
the third statement holds, then the set of principal
$G$--bundles $\{(F^i_M)^*E_G\}_{i\geq 1}$ contains
only finitely many isomorphism classes. Hence the
first statement in the theorem follows from the third
statement.

Now assume that the second statement in the theorem holds.
Therefore, for each $G$--module $V$, the vector bundle
$E_G(V)\, \longrightarrow\, M$
associated to $E_G$ for $V$ is essentially finite.
We note that for an essentially finite vector bundle
$W\, \longrightarrow\, M$, there is a finite group--scheme
$\Gamma$ over $k$ and a principal $\Gamma$--bundle
$$
\gamma\, :\, E_\Gamma\, \longrightarrow\, M
$$
such that the vector bundle $\gamma^*W
\, \longrightarrow\, E_\Gamma$ is trivializable;
this follows from \cite[p. 83, Proposition 3.10]{No}. Using
this we conclude the following:
\begin{itemize}
\item $W$ is strongly semistable, and

\item{} $[c_i(W)c_1(\xi)^{d-i}]\, =\, 0$
for all $i\, \geq\, 1$.
\end{itemize}
In particular,  $E_G(V)$ is strongly semistable and
$$
[c_i(E_G(V))c_1(\xi)^{d-i}]\, =\, 0
$$
for all $i\, \geq\, 1$. Consequently, the third
statement in the theorem holds.

The assertion that
the third statement implies the second statement is
essentially contained in
Lemma \ref{lem1}. To explain this, assume that
the third statement holds.

Since the field $k$ is finite, a theorem of Lang says that
the fiber $(E_G)_{x_0}$ is a trivial $G$--torsor (see
\cite[p. 557, Theorem 2]{Lang}). Fix a $k$--rational point
\begin{equation}\label{e6}
z_0\, \in\, (E_G)_{x_0}\, .
\end{equation}

Let $\text{Rep}(G)$ denote the category of all finite
dimensional representations of $G$ over $k$. For any $V\, \in\,
\text{Rep}(G)$, let $E_G(V)$ be the
vector bundle over $M$ associated to $E_G$ for $V$.
The point $z_0$ in Eq. \eqref{e6} defines an isomorphism
\begin{equation}\label{e7}
f_{V,z_0}\, :\, V\, \longrightarrow\, E_G(V)_{x_0}
\end{equation}
of vector spaces over $k$. More precisely, $f_{V,z_0}$ sends
any vector $w\, \in\, V$ to the image of $(z_0\, ,w)$ in
$E_G(V)_{x_0}$ (recall that $E_G(V)\, :=\, E_G\times^G V$
is a quotient of $E_G\times V$).

Let
\begin{equation}\label{e8}
{\mathcal F}_0\, :\, \text{Rep}(G)\,\longrightarrow\,
\text{Rep}(\varpi(M,x_0))
\end{equation}
(see Eq. \eqref{e5}) be the functor that sends any $G$--module
$V$ to the essentially finite vector bundle $E_G(V)$ associated
to $E_G$ for $V$ (it was shown in Lemma \ref{lem1} that
$E_G(V)$ is essentially finite); recall that an
essentially finite vector bundle over $M$ gives an object of
$\text{Rep}(\varpi(M,x_0))$.

Let $k$--mod denote the category of finite dimensional vector
spaces over $k$. Let
$$
{\mathcal G}\, :\, \text{Rep}(G)\,\longrightarrow\,
k\mbox{--}{\rm mod}
$$
be the fiber functor for $G$ that sends any $G$--module to the
underlying $k$--vector space. Similarly, let
$$
{\mathcal H}\, :\, \text{Rep}(\varpi(M,x_0))\,\longrightarrow\, 
k\mbox{--}{\rm mod}
$$
be the fiber functor for $\varpi(M,x_0)$ that sends
an essentially finite vector bundle $W$ to its fiber
$W_{x_0}$ over $x_0$. For any $V\, \in\,
\text{Rep}(G)$, the homomorphism $f_{V,z_0}$ in Eq. \eqref{e7}
is an isomorphism of ${\mathcal G}(V)$ with
${\mathcal H}({\mathcal F}_0(V))$, where ${\mathcal F}_0$ is
constructed in Eq. \eqref{e8}.

Now we note that ${\mathcal F}_0$ and the
isomorphisms $\{f_{V,z_0}\}_{V\in
\text{Rep}(G)}$ together define a homomorphism
of group--schemes
\begin{equation}\label{rho}
\rho\,:\, \varpi(M,x_0)\, \longrightarrow\, G
\end{equation}
(recall that $\varpi(M,x_0)$ is defined to be the group--scheme
corresponding to the neutral Tannakian category defined
by the category of essentially finite vector bundles on $M$
equipped with the fiber functor that sends any 
essentially finite vector bundle to its fiber over $x_0$).

Let $E_\rho$ denote the principal $G$--bundle over $M$ obtained
by extending the structure group of the tautological principal
$\varpi(M,x_0)$--bundle $E_{\varpi(M,x_0)}\, \longrightarrow\, M$
using
the homomorphism $\rho$ in Eq. \eqref{rho}. Note that the morphism
$$
\text{Id}_{E_{\varpi(M,x_0)}}\times \rho\, :\,
E_{\varpi(M,x_0)}\times \varpi(M,x_0)\, \longrightarrow\,
E_{\varpi(M,x_0)}\times G
$$
descends to a morphism
\begin{equation}\label{ad.ex.h} 
\widetilde{\rho}\, :\,
E_{\varpi(M,x_0)}\, =\, E_{\varpi(M,x_0)}\times^{\varpi(M,x_0)}
\varpi(M,x_0)\, \longrightarrow\, E_{\varpi(M,x_0)}\times^{\varpi(M,x_0)}
G\, =:\, E_\rho
\end{equation}
between the quotient spaces.

{}From the construction of the homomorphism $\rho$ in
Eq. \eqref{rho} it follows that for any $V\, \in\,\text{Rep}(G)$,
the vector bundle $E_\rho\times^G V\, \longrightarrow\,
M$ associated to $E_\rho$ for $V$
is identified with the vector bundle $E_G(V)$ associated
to $E_G$ for $V$. Therefore, we get an isomorphism of
the principal $G$--bundle $E_\rho$
with $E_G$. This completes the proof of the theorem.
\end{proof}

Note that the above identification of $E_G$ with $E_\rho$
takes the point in $E_\rho$ defined by $(e_0\, ,e)\,
\in\, E_{\varpi(M,x_0)}\times G$, where $e\, \in\, G$
is the identity element and $e_0$ is the element in Eq. \eqref{e4},
to the point $z_0$ of $E_G$ in Eq. \eqref{e6}.

The homomorphism $\rho$ in Eq. \eqref{rho} depends on the choice
of $z_0$. Take any $g\, \in\, G$. Let 
$$
\rho'\,:\, \varpi(M,x_0)\, \longrightarrow\, G
$$
denote the homomorphism obtained in place of $\rho$ after
we replace $z_0$ by $z_0'\, =\, z_0g$. Then,
\begin{equation}\label{e9}
\rho'\, =\, I_{g_0}\circ \rho\, ,
\end{equation}
where $I_{g_0}$ is the inner automorphism of $G$ defined by
$g\,\longmapsto\, g^{-1}_0gg_0$.

Consider the adjoint action of $G$ on itself.
Let
$$
\text{Ad}(E_G)\, =\, E_G\times^G G\, \longrightarrow\,
M
$$
be the associated fiber bundle. Since the adjoint action
of $G$ on itself preserve the group structure, it follows
that $\text{Ad}(E_G)$ is a group--scheme over $M$.
Let
$$
\text{Ad}(E_{\varpi(M,x_0)})\, :=\, E_{\varpi(M,x_0)}
\times^{\varpi(M,x_0)} \varpi(M,x_0)\, \longrightarrow\, M
$$
be the adjoint group--scheme
for $E_{\varpi(M,x_0)}$. Since
$E_G\,=\, E_\rho$ is an extension of structure group of
$E_{\varpi(M,x_0)}$ using $\rho$, we have a homomorphism
\begin{equation}\label{e10}
\varphi\, :\, \text{Ad}(E_{\varpi(M,x_0)})
\, \longrightarrow\, \text{Ad}(E_G)
\end{equation}
of group--schemes over $M$. Indeed, the morphism
$$
\widetilde{\rho}\times \rho\, :\, E_{\varpi(M,x_0)}\times
\varpi(M,x_0)\,\longrightarrow\, E_\rho\times G\, ,
$$
where $\widetilde{\rho}$ is constructed in Eq. \eqref{ad.ex.h},
descends to the morphism $\varphi$ in Eq. \eqref{e10} between 
the quotient spaces. Now, using Eq. \eqref{e9} it is
straight forward to check that the morphism
$\varphi$ does not depend on the choice of the point $z_0$. The
image of $\varphi$ is also independent of the choice
of the base point $x_0$.

Lemma \ref{lem1} and Theorem \ref{thm1} have the following
corollary:

\begin{corollary}\label{cor1}
A vector bundle $E$ over $M$ is essentially finite if and only
if there are integers $b\, >\, a\, \geq\, 1$, such that
the vector bundle $(F^b_M)^*E$ is isomorphic to $(F^a_M)^*E$.
\end{corollary}

\begin{proof}
Assume that $(F^b_M)^*E$ is isomorphic to $(F^a_M)^*E$ for some
$a$ and $b$ with $b\, >\, a\, \geq\, 1$. Hence Eq. \eqref{eqf}
holds. Now exactly as in the proof of Lemma \ref{lem1}, using the 
theorem
of Lange--Stuhler and Deligne together with Proposition 2.3 of 
\cite{BH} we conclude that $E$ is essentially finite.

Now assume that $E$ is an essentially finite vector bundle of
rank $r$. Set $G\, =\, \text{GL}(r,k)$, and substitute for
$E_G$ in Theorem \ref{thm1}
the principal $\text{GL}(r,k)$--bundle
$E_{\text{GL}(r,k)}\, \longrightarrow\, M$ given by $E$. Since
$E$ is essentially finite, we know that there is homomorphism
$$
\rho\, :\, \varpi(M,x_0)\, \longrightarrow\, \text{GL}(r,k)
$$
such that $E_{\text{GL}(r,k)}$ is isomorphic to the principal
$\text{GL}(r,k)$--bundle
obtained by extending the structure group of the
tautological principal $\varpi(M,x_0)$--bundle
$E_{\varpi(M,x_0)}\, \longrightarrow\, M$ using $\rho$;
the homomorphism $\rho$ is given by the functor
$$
\text{Rep}(\text{GL}(r,k))\, \longrightarrow\,
\text{Rep}(\varpi(M,x_0))
$$
that sends a $\text{GL}(r,k)$--module $V$ to the associated
vector bundle $E_{\text{GL}(r,k)}(V)$.
Hence from Theorem \ref{thm1} we know that there
are integers $b\, >\, a\, \geq\, 1$, such that the principal
$\text{GL}(r,k)$--bundle $(F^b_M)^*E_{\text{GL}(r,k)}$
is isomorphic to $(F^a_M)^*E_{\text{GL}(r,k)}$.
This completes the proof of the corollary.
\end{proof}



\begin{thebibliography}{1111}

\bibitem{BPS} Biswas, I., Parameswaran, A. J. and
Subramanian, S.: Monodromy group
for a strongly semistable principal bundle over a curve,
\textit{Duke Math. Jour.} \textbf{132} (2006), 1--48.

\bibitem{BH} Biswas, I. and Holla, Y. I.: Comparison
of fundamental
group schemes of a projective variety and an ample hypersurface,
\textit{Jour. Alg. Geom.} \textbf{16} (2007), 547--597.

\bibitem{Bo} Bogomolov, F. A.: Holomorphic tensors and vector
bundles on projective varieties, \textit{Math. USSR--Izv.}
\textbf{13} (1978), 495--555.

\bibitem{HK} Coiai, F. and Holla, Y. I.: Extensions of structure
groups of principal bundles in positive characteristics,
\textit{Jour. Reine Angew. Math.} \textbf{595} (2006), 1--24.

\bibitem{DM} Deligne, P. and Milne, J. S.: Tannakian Categories,
in: \textit{Hodge cycles, motives, and Shimura varieties} (by P.
Deligne, J. S. Milne, A. Ogus and K.-Y. Shih), pp. 101--228,
Lecture Notes in Mathematics, \textbf{900}, Springer-Verlag,
Berlin-Heidelberg-New York, 1982.

\bibitem{Gi} Giraud, J.: Cohomologie non ab\'elienne, Die Grundlehren
der mathematischen Wissenschaften, Band 179, Springer--Verlag,
Berlin--New York, 1971. 

\bibitem{Lang} Lang, S.: Algebraic groups over finite fields,
\textit{Amer. Jour. Math.} \textbf{78} (1956), 555--563.

\bibitem{LS} Lange, H. and Stuhler, U.: Vektorb\"undel auf 
Kurven und Darstellungen der algebraischen Fundamentalgruppe,
\textit{Math. Zeit.} \textbf{156} (1977), 73--83.

\bibitem{Lan} Langer, A.: Semistable sheaves in positive
characteristic, \textit{Ann. of Math.} \textbf{159}
(2004), 251--276.

\bibitem{Lan2} Langer, A.: Semistable principal $G$--bundles
in positive characteristics, \textit{Duke Math. Jour.}
\textbf{128} (2005), 511--540.

\bibitem{La} Laszlo, Y.: A non--trivial family of bundles
fixed by the square of Frobenius, \textit{Comp. Ran. Acad.
Sci. Paris} \textbf{333} (2001), 651--656.

\bibitem{No0} Nori, M. V.: On the representations of
the fundamental group scheme, \textit{Compos. Math.} 
\textbf{33} (1976), 29--41

\bibitem{No} Nori, M. V.: The fundamental
group--scheme, \textit{Proc. Ind. Acad. Sci. (Math. Sci.)}
\textbf{91} (1982), 73--122.

\bibitem{RR} Ramanan, S. and Ramanathan, A.: Some
remarks on the instability flag, \textit{T\^ohoku Math. Jour.}
\textbf{36} (1984), 269--291.

\bibitem{Ra} Ramanathan, A.: Stable principal bundles on a
compact Riemann surface, \textit{Math. Ann.}
\textbf{213} (1975), 129--152.

\bibitem{Su} Subramanian, S.: Strongly semistable bundles on
a curve over a finite field, \textit{Arch. Math.}
\textbf{89} (2007), 68--72.

\end{thebibliography}
\end{document}